\documentclass[a4paper, 12pt]{article}
\usepackage{amsfonts}
\usepackage{mathrsfs}
\usepackage{amsmath}

\usepackage{amsthm}
\usepackage{array}
\usepackage{cases}
\makeatletter
\def\th@plain{%
  \upshape 
}
\makeatother



\newtheorem{thm}{Theorem}

\newtheorem{lem}[thm]{Lemma}

\numberwithin{equation}{section}

\usepackage[varg]{txfonts}
\usepackage{graphicx, epsfig, subfigure}
\usepackage{enumerate} 
\usepackage[square, numbers, sort&compress]{natbib} 
\ifx\pdfoutput\undefined
 \usepackage[dvipdfm,%
  pdfstartview=FitH, 
  bookmarks=true,%
  bookmarksnumbered=true, 
  bookmarksopen=true, 
  plainpages=false,%
  pdfpagelabels,%
  colorlinks=true, 
  linkcolor=blue, 
  citecolor=blue,%
  urlcolor=black,
  pdfborder=001]{hyperref}
  \AtBeginDvi{}  
\else
 \usepackage[pdftex,%
  pdfstartview=FitH, 
  bookmarks=true,%
  bookmarksnumbered=true, 
  bookmarksopen=true, 
  plainpages=false,%
  pdfpagelabels,%
  colorlinks=true, 
  linkcolor=blue, 
  citecolor=blue,%
  urlcolor=black,
  pdfborder=001]{hyperref}
\fi

\numberwithin{equation}{section}

\begin{document}
\title{Class two 1-planar graphs with maximum degree six or seven}
\author{Xin Zhang\thanks{Email addresses: sdu.zhang@yahoo.com.cn, xinzhang@mail.sdu.edu.cn.}~\thanks{Supported by GIIFSDU (yzc10040) and NSFC (10971121).}\\[.5em]
{\small School of Mathematics, Shandong University, Jinan 250100, P. R. China}\\
}
\date{}
\maketitle

\begin{abstract}\baselineskip 0.60cm
A graph is $1$-planar if it can be drawn on the plane so that each edge is crossed by at most one other edge. In this note we give examples of class two 1-planar graphs with maximum degree six or seven.
\\[.5em]

\end{abstract}
\baselineskip 0.60cm

A graph is 1-planar if it can be drawn on the plane so that each edge is crossed by at most one other edge. This notion of $1$-planar graphs was introduced by Ringel \cite{Ringel.1965} while trying to simultaneously color the vertices and faces of a planar graph $G$ such that any pair of adjacent/incident elements receive different colors. The coloring problems of 1-planar graphs have been investigated in many papers such as \cite{Albertson.2006,Borodin.1984,Borodin.1995,Borodin.2001,Zhang.2010.SDU,Zhang.2011,Zhang.new}.

Vizing's theorem (see page 251 of \cite{Book}) states that the edge chromatic number of every nonempty graph $G$ is either $\Delta(G)$ or $\Delta(G)+1$. Thus we can divide all graphs into two classes. A graph $G$ is of Class one if $\chi'(G)=\Delta(G)$ and is of Class two if $\chi'(G)=\Delta(G)+1$. Consequently, a major question in the area of edge colorings is that of determining to which of these two classes a given
graph belongs.

For a 1-planar graph $G$, it is proved that $G$ is of Class one provided $\Delta(G)\geq 10$ \cite{Zhang.2011}, or $\Delta(G)\geq 8$ and $G$ contains no adjacent triangles \cite{Zhang.new}, or $\Delta(G)\geq 7$ and $G$ is triangle-free \cite{Zhang.2010.SDU}. Moreover, in \cite{Zhang.new}, the authors conjectured that every 1-planar graph with maximum degree at least 8 is of Class 1.

In this note we aim to construct examples of class two 1-planar graphs with maximum degree six or seven.
Note that every planar graph is 1-planar graph and Vizing \cite{Vizing} presented examples of planar graphs of Class two with maximum degree no more than five. Therefore, we can conclude that
there are 1-planar graphs of Class two with maximum degree $\Delta$ for each $\Delta\leq 7$.

At first, we need an useful lemma. Recall that $\alpha'(G)$ is the edge independent number of $G$ and $\alpha'(G)\leq \lfloor\frac{n}{2}\rfloor$, where $n$ is the order of $G$.

\begin{lem}(see page 258 of \cite{Book})\label{lem}
If $G$ is a graph of size $m$ such that $m>\alpha'(G)\Delta(G)$, then $G$ is of Class two.
\end{lem}

Now we check that the above two 1-planar graphs are of Class two by Lemma \ref{lem}. It is easy to see that the order of $G_1$ is 25 and the size of $G_1$ is 73. Thus $\alpha'(G_1)\leq 12$, which implies that $|E(G_1)|=73>72\geq \alpha'(G_1)\Delta(G_1)$. Similarly, one can also see that $|E(G_2)|=85>84=\lfloor\frac{|V(G_2)|}{2}\rfloor\Delta(G_2)\geq \alpha'(G_2)\Delta(G_2)$. So by Lemma \ref{lem}, $G_1$ and $G_2$ are both of Class two.

Note that all vertices besides one 2-vertex in $G_1$ (resp$.$ $G_2$) are of degree 6 (resp$.$ 7). So we can respectively construct a 6-regular and a 7-regular 1-planar graph by deleting the unique 2-vertex and adding an edge between the two neighbors of it. Here we would like to point out that the resulting 7-regular 1-planar graph was presented by Fabrici and Madaras in their paper \cite{Fabrici.2007}.

\begin{figure}
\begin{center}
  \includegraphics[width=6cm,height=5.5cm]{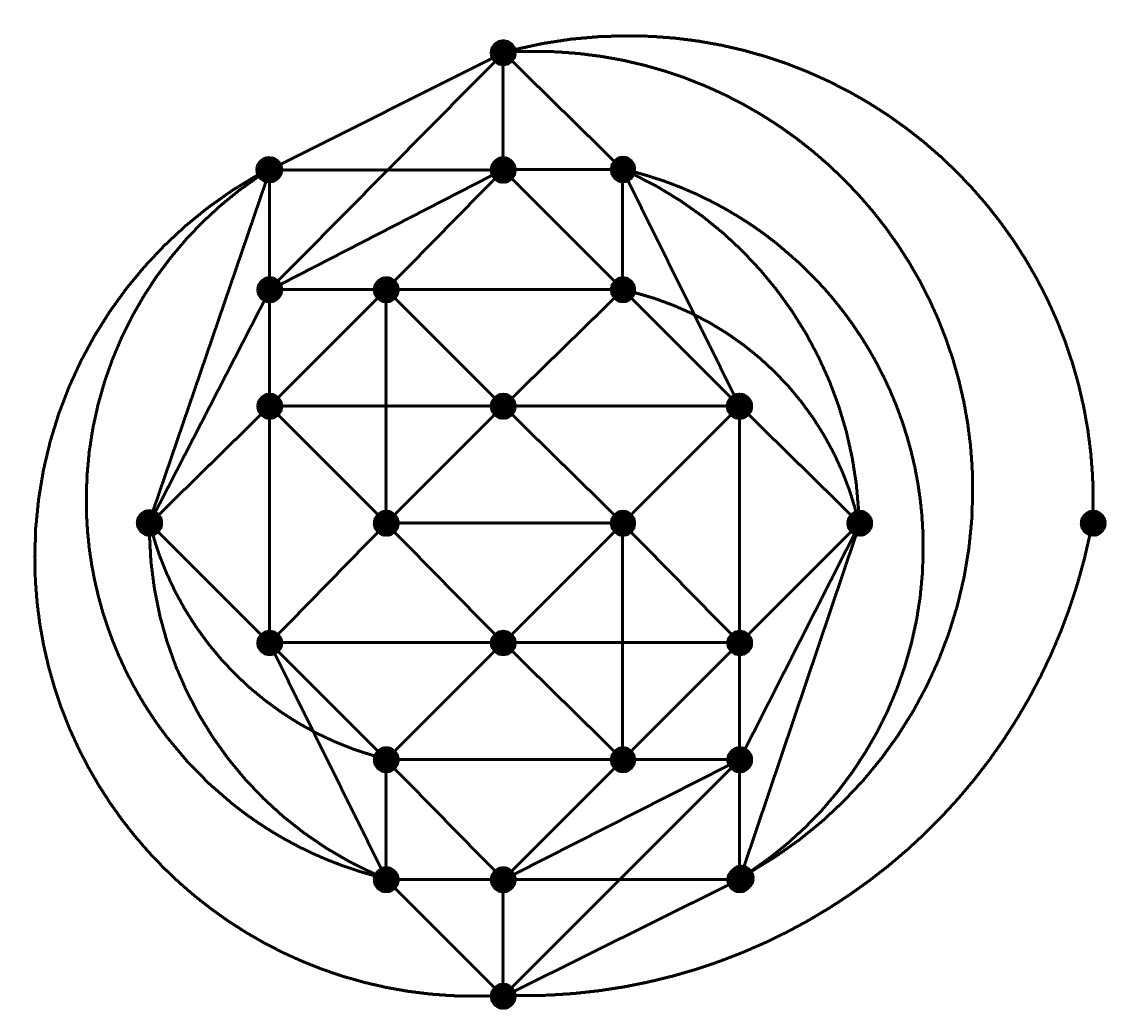}~~~~~
  \includegraphics[width=6cm,height=5.5cm]{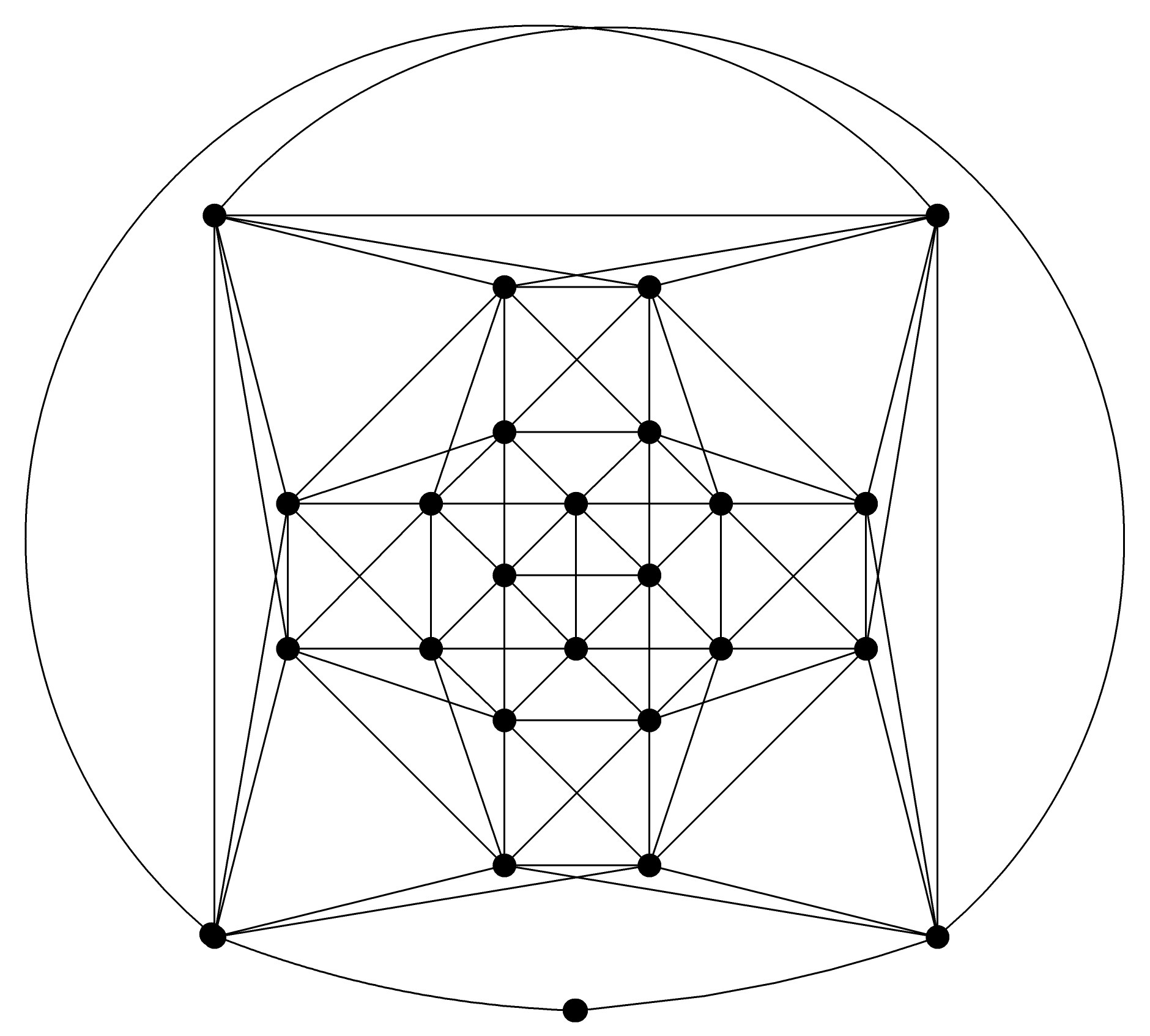}\\
  \caption{Class two 1-planar graphs $G_1$ with maximum degree 6 and $G_2$ with maximum degree 7.}\label{class2}
\end{center}
\end{figure}


\begin{thebibliography}{10}\setlength{\itemsep}{0pt}


\bibitem{Albertson.2006} M. O. Albertson, B. Mohar. Coloring vertices and faces of locally planar graphs. Graphs and Combinatorics \textbf{22} (2006) 289--295.

\bibitem{Borodin.1984} O. V. Borodin. Solution of Ringel's problems on the vertex-face coloring of plane graphs and on the coloring of $1$-planar graphs. Diskret. Analiz \textbf{41} (1984) 12--26.

\bibitem{Borodin.1995} O. V. Borodin. A New Proof of the $6$-Color Theorem. Journal of Graph Theory \textbf{19(4)} (1995) 507--521.

\bibitem{Borodin.2001} O. V. Borodin, A. V. Kostochka, A. Raspaud, E. Sopena. Acyclic colouring of 1-planar graphs. Discrete Applied Mathematics \textbf{114} (2001) 29--41.

\bibitem{Book} G. Chartrand, P. Zhang. Chromatic Graph Theory. CRC Press, 2008.

\bibitem{Fabrici.2007} I. Fabrici, T. Madaras. The structure of 1-planar graphs. Discrete Mathematics \textbf{307} (2007) 854--865.

\bibitem{Ringel.1965} G. Ringel. Ein Sechsfarbenproblem auf der Kugel. Abh. Math. Sem. Univ, Hamburg \textbf{29} (1965) 107--117.

\bibitem{Vizing} V. G. Vizing, Critical graphs with given chromatic class, Diskret. Analiz. \textbf{5} (1965) 9--17.

\bibitem{Zhang.2010.SDU} X. Zhang, G. Liu, J.-L. Wu. Edge coloring of triangle-free 1-planar graphs. Journal of Shandong
University (Natural Science) \textbf{45(6)} (2010) 15--17.

\bibitem{Zhang.2011} X. Zhang, J.-L. Wu. On edge colorings of 1-planar graphs. Information Processing Letters \textbf{111(3)} (2011) 124--128.

\bibitem{Zhang.new} X. Zhang, G. Liu. On edge colorings of 1-planar graphs without adjacent triangles. Submitted.

\end{thebibliography}
\end{document}